\RequirePackage{fix-cm}
\documentclass[smallextended]{svjour3}       
\smartqed  
\usepackage{graphicx}

\usepackage{amsmath,amssymb,amsbsy}

\usepackage{latexsym}
\begin{document}

\title{The time-adaptive    statistical testing for random number generators
}


\author{Boris Ryabko
}


\institute{ \at  Institute of  Computational Technologies  of
Siberian Branch of Russian Academy of Science \\
and Novosibirsk State University. 
              \email{boris@ryabko.net}           
         }

\date{Received: date / Accepted: date}

\maketitle

\begin{abstract}
The problem of constructing effective statistical tests for random number generators (RNG) is considered. Currently, there are hundreds of RNG statistical tests that are often combined into so-called batteries, each containing from a dozen to more than one hundred tests.  
 When a battery test is used,  it  is applied to a sequence generated by the RNG, and the calculation time is determined by the length of the sequence and the number of tests. Generally speaking, the longer the sequence, the smaller deviations from randomness can be found by a specific test. So, when a battery is applied, on the one hand, the ``better'' tests are  in the battery, the more chances to reject a ``bad'' RNG. On the other hand, the larger the battery, the less time can be spent on each test and, therefore, the shorter the test sequence. In turn, this reduces the ability to find small deviations from randomness.
To reduce this trade-off, we propose  an adaptive way to use batteries (and other sets) of tests, which requires less time but, in a certain sense, preserves the power of the original battery. We call this method  time-adaptive   battery of tests. 

The suggested method is 
based on the theorem which describes asymptotic  properties of the so-called p-values of tests. Namely, the theorem claims that if the RNG can be modelled by a stationary ergodic source, the  value $- log \, \pi (x_1 x_2 ... x_n)  /n$ goes to 
$1 - h$ when $n$ grows,  where $x_1 x_2 ... $ is the  sequence, $\pi (\, )$ is the p- value of the most powerful  test, $h$ is the  limit Shannon entropy of the stationary ergodic source.

\end{abstract}

\section{Introduction}  
Random number generators (RNG) and pseudo-random number generators (PRNG) are widely used in many applications. 
   RNGs are based on physical sources, while pseudo-random numbers are  generated by computers.
The goal of RNG and PRNG is to generate sequences of binary digits, which are distributed as a result of throwing an ``honest'' coin, or, more precisely, obey the Bernoulli distribution with parameters $ (1/2, 1/2) $. As a rule, for practically used RNG and PRNG this property is verified experimentally with the help of statistical tests developed for this purpose. 
 
Currently, there are more than one hundred applicable 
statistical tests, as well as  dozens RNGs based on different physical processes, and an even greater number of PRNGs based on different mathematical algorithms;  see for review \cite{l2017history,l2007testu01,grev}. 
Informally, an ideal RNG should generate sequences that pass all tests. In practice, especially in cryptographic applications, this requirement is formulated as follows: an RNG must pass a so-called battery of statistical tests, that is, some fixed set of tests. 
When a battery is applied,  each test in the
test battery is applied separately to the RNG. 
Among these batteries, we mention the Marsaglia's Diehard battery, which contains 16 tests \cite{mar}, the National Institute of Standards and Technology (NIST) battery of 15 tests \cite{NIST-test}, several   batteries proposed by L'Ecuyer and Simard \cite{l2007testu01}, which contain from 10 to 106 tests
 and  many  others (see for review \cite{l2017history,l2007testu01,tur}).  In addition, these batteries contain many tests that can be used with different  values of the parameters, potentially increasing the total number of tests in the battery.  Note that practically used RNG should be tested from time to time like any physical equipment, and therefore these test batteries should be used continuously.

How to evaluate  large  batteries of tests? On the one hand, the larger the test battery, the more likely it is to find flaws in the tested RNG. On the other hand, the larger the battery, the more time is required for testing. 
(Thus, L'Ecuyer and Simard \cite{l2007testu01}
remark the need for
small batteries to increase computational efficiency.)
   Another view is 
as follows: in reality,  the  time available  to study any RNG is limited.  Given a  certain time budget, 
 one  can either   use more tests and  relatively short sequences generated by the RNG, or use fewer tests,  but longer sequences and, in turn, this gives more chances to find deviations of randomness of the considered RNG. 

In order to  reduce this trade-off, we propose a time-adaptive  testing of RNGs, in which, informally speaking, first all the tests are executed on  relatively short sequences generated by the RNG, and then a few ``promising'' tests are applied for the final testing.
Of course, the key question here is which tests are promising. For example, if a battery of two tests is applied to (relatively short) sequences of the same length, it can be assumed that the smaller  the p-value, the more promising the test. But a more complicated situation may arise when we have to compare two tests that were applied to sequences of different lengths (for example, the first test was applied to a sequence of length $ l_1 $, and the second to a sequence  of length of $ l_2 $, $ l_1 \neq l_2 $). We show that if our goal is to choose the most powerful test, then a good strategy is to choose the test $i$  for which the ratio $- \log (p-value_i) / l_i $ is maximum.
This 
recommendation is based on the following theorem:
 if an RNG can be modelled by a stationary ergodic source, the  value $- log \, \pi (x_1 x_2 ... x_n)  /n$ goes to 
$1 - h$, if $n$ grows,  where $x_1 x_2 ... $ is a generated sequence, $\pi (\, )$ is the p-value of the most powerful  test, $h$ is the  limit Shannon entropy of the stationary ergodic source.   This theorem plays an important rule in the suggested time-adaptive scheme and will be described in the first part of the paper, whereas  the time-adaptive testing will be described afterwords. The description will be illustrated by experiments with the battery Rabbit  from \cite{l2007testu01}.

As far as we know, the proposed approach to testing RNGs is new, but the idea of finding the best test among many, testing the tests step by step in an increasing sequence, is widely used in  algorithmic information theory, where the notion of random sequence is formally investigated and discussed \cite{Calude:02,downey2006calibrating}. 


\section{Hypothesis testing and properties of pi-values}

\subsection{Notation}

We consider RNG which generates a sequence of  letters $x = x_1 x_2 $ $ ... x_n, $ $n \ge 1,$  from a finite alphabet  $ \{0,1\}^n$. Two  statistical hypotheses are considered: $H_0 = $ 
$\{ x \, $ \emph{ 
obeys the uniform distribution  ($ \mu_U$) on  
  $\{0,1\}^n$ \}},  and the alternative hypothesis $H_1 = \bar{H_0} 	$, 
that is, $H_1$ is the negation of $H_0$.   It is a particular case of the so-called  goodness-of-fit problem, and any test for it   is called a test of fit, see \cite{ks}.
Let $t$ be a test. Then,  by definition, the significance level $\alpha$ equals the probability 
of the Type I error, $\alpha \in (0,1)$.  
Denote a critical region of the test $t$ for the  significance level $\alpha$ by $C_t(\alpha)$  and let 
$\bar{C}_t(\alpha)$ $ =  \{0,1\}^n \setminus C_t(\alpha) .$ 
(Recall that Type I error occurs   if
$H_0$ is true and is rejected. Type II error occurs   if
$H_1$ is true, but $H_0$ is accepted.   Besides, 
for a certain 
$x = x_1 x_2 ... x_n\,\,$
$H_0$ is rejected if and only if $x \in C_t(\alpha)$.)

Suppose that $H_1$ is true, and  the investigated sequence 
 $x = x_1 x_2 ... x_n$ is generated by an 
(unknown)   source $\nu$.
By definition, a test $t$ is consistent, if for any 
significance level $\alpha \in (0,1)$ the probability of Type II error goes to 0, that is 
 \begin{equation}\label{cons}
 \lim_{n \to \infty} \nu( \bar{C}_t(\alpha) ) = 0 \, .
 \end{equation}

Suppose, that $H_1$ is true and the sequences $x \in \{0,1\}^n$  obey a certain distribution $\nu$.
It is well-known in mathematical statistics that  the optimal test  (Neyman-Pearson or $NP$ test)  is described by the
Neyman-Pearson lemma and the critical region of this test  is defined 
as follows:
$$
C_{NP}(\alpha) = \{x :  \, \, \mu_U (x) / \nu (x) \le \lambda_\alpha  \} \, ,
$$
where $\alpha \in (0,1) $ is the 
significance level and the constant $\lambda_\alpha$ is chosen in such a way that 
$\mu_U ( C_{NP}(\alpha) ) = \alpha$, see \cite{ks}.
(We did not take into account that the set $ \{0,1\}^n$ is finite. 
Strictly speaking, in such a case a randomized test should be used, but in what follows 
we will consider asymptotic behaviour of  tests for large $n$, and this effect will be negligible).
Note that, by definition, $ \mu_U (x)  = 2^{-n}$ for any $x$ $\in \{0,1\}^n$.

\subsection{The p-value and its properties.}
The notion of the critical region is connected with the so-called p-value, which we define for the NP-test by the following equation:
 \begin{equation}\label{pi-np}
\pi_{NP} (x) = \mu_U \{ y: \nu(y) > \nu(x) \} = | \{ y: \nu(y) > \nu(x) \} | / 2^n \, .
\end{equation}
Informally, 
$ \pi_{NP} (x)$ is the probability to meet a random point $y$ which is 
worse than the
observed when considering the null hypothesis.

 The NP-test is optimal in the sense that its probability of a Type II error is minimal, but when testing an RNG the alternative distribution is unknown, and, hence,   different tests are necessary. 
 Let us consider a certain  statistic $\tau$  (that is, a function on $\{0,1\}^n$),
and define the  p-value for this $\tau$ and  $x$ as follows:
 \begin{equation}\label{pi-tau}
\pi_\tau (x) = \mu_U \{ y: \tau(y) > \tau(x) \} = | \{ y: \tau(y) > \tau(x) \} | / 2^n \, .
 \end{equation}
(Note, that  the definition $\pi_{NP} $ in  (\ref{pi-np}) corresponds to this equation if
the value $\nu(x)$ is considered as a statistic, i.e. $\tau(x) = \nu(x)$).

 \subsection{The p-value and Shannon entropy.}  
 
 It turns out that there exist such  tests whose asymptotic behaviour is close to that of the $NP$-test 
 for any (unknown) stationary ergodic source $\nu$, see \cite{arx}. 
Those tests are based on so-called universal codes (or data-compressors)  and 
 are described in \cite{R1,R2},  where it is shown that  they are consistent. 
We  describe those tests  in Appendix 1 and show that they are asymptotically optimal.
 The following theorem  describes the asymptotic behaviour of p-values for stationary ergodic sources for $NP$ test and the mentioned above tests which are based on universal codes (see Appendix 1). We use this theorem as the theoretical  basis for adaptive statistical testing developed in this paper.  
\begin{theorem}\label{t1} 
i)  If $\nu$ is a stationary ergodic measure, then, with probability 1, 
\begin{equation}\label{t11}
 \lim_{n \to \infty} -  \frac{1}{n} \log \pi_{NP} (x) = 1 - h(\nu) \, ,
\end{equation}
where $h(\nu)$ is the Shannon entropy of $\nu$, see for definition \cite{co}.

ii) There exists such  a statistic $\tau$  that for any  stationary ergodic measure  
$\nu$, with probability 1, 
\begin{equation}\label{t2}
 \lim_{n \to \infty} -  \frac{1}{n} \log \pi_{\tau} (x) = 1 - h(\nu) \, ,
\end{equation}
where p-values $ \pi_{NP}$ and $\pi_{\tau}$ are  defined in (\ref{pi-np}) and  (\ref{pi-tau}), 
correspondingly.
\end{theorem}

The statistic $\tau$  and the corresponding test of fit 
are described in  Appendix 1,  
the proof of the theorem   is given in Appendix 2, 
but here we   note that this theorem gives some idea of the relation between the Shannon entropy of the (unknown) process $ \nu $ and the required sample size.
Indeed,  suppose that the $ NP $ test is used and the desired significance level is $ \alpha $. Then, we can see  that
(asymptotically) $\alpha$ should be larger than $\pi_{NP}(x)$ and from (\ref{t11}) 
we obtain 
$ n > - \log \alpha / (1-h(\nu) ) \, $ (for the most powerful  test).
It is known that the Shannon entropy is 1 if and only if $\nu$ is a uniform measure $\mu_u$.  Therefore, in a certain sense, the difference $ 1-h (\nu) $ estimates the distance between the distributions, and the last inequality shows that the sample size becomes infinite if $ \nu $ approaches a uniform distribution.

The next simple example illustrates  this theorem.
 Let there be a statistic $\tau$ and a generator (a measure $\nu$)  created sequences of binary digits which 
are independent and, say,   $\nu(0) = 0.501, \nu(1) = 0.499 $.
Suppose, 
 $\lim_{n \to \infty} -  \frac{1}{n}  $ $ \log \pi_{\tau} (x) $ $ = c \, $,
 where $c$ is a positive constant.  Let us consider the following  ``decimation test''  $\tau^{1/2}$:  an input sequence $x_1 x_2 .... x_n$ is transformed into $x_1 x_3 x_5 ... x_{2 \lfloor n/2 \rfloor -1} $ and then  the test is applied to this transformed sequence. Obviously, for this test 
 $\lim_{n \to \infty} -  \frac{1}{n/2} \log \pi_{\tau^{1/2}} (x) = c  \, $, and, 
 hence, $\lim_{n \to \infty} $ $-  \frac{1}{n} $ $ \log \pi_{\tau^{1/2}} (x) = c /2 \, $. 
 Thus, the value $ - \frac {1} {n} \log \pi_{\tau} (x_1 ... x_n) $ seems to be a reasonable estimate of the power of the test   for a large $n$.
 

\section{Time-adaptive statistical tests and their experimental investigation }
\subsection{Batteries of tests.} 

Let us consider a situation where the randomness testing  is performed by conducting  a battery  of statistical  tests for randomness.  Suppose that the battery contains $s$ tests  and  $\alpha_i$ is the significance level of $i-$th test, $i= 1, ... , s$.  If the battery is applied in such a way that
 the hypothesis $H_0$ is rejected when at least one  test in the battery rejects it, then the significance level $\alpha$
of this battery satisfies the following inequality:
 \begin{equation}\label{a-g}
\alpha \le \sum_{i=1}^s \alpha_i \, .
\end{equation}
If all the tests in the battery are independent, then the following  equation is valid:
$
\alpha = 1 - \prod_{i=1}^s (1- \alpha_i) \, .
$  Clearly,  the upper bound (\ref{a-g})  is true  for this case and $1 - \prod_{i=1}^s (1- \alpha_i)$ is close to $\sum_{i=1}^s \alpha_i$, if each $\alpha_i$ is much smaller
than $1/s$. That is why we will use the estimate (\ref{a-g}) below.

We have considered a scenario in which a test  is applied to a single sequence generated by an RNG, and then the researcher makes a decision on the RNG based on the test results.  Another possibility that has been considered
by several authors, e.g.\ \cite{l2007testu01,NIST-test}, 
is to use  the following two-step 
procedure for testing RNGs. 
The idea is to generate $ r $ sequences $ x ^ 1, x ^ 2, ..., x ^ r $ and apply one test (say, $ \tau $) to each of them independently.
Then apply another  test to the received data $ \tau (x ^ 1), \tau (x ^ 2), ..., \tau (x ^ r) $ (as a rule,  those 
values are converted  into a sequence  of corresponding p-values, and then the hypothesis of the uniform distribution of those p-values   is tested).
  Then this procedure is repeated for the second test in the battery, and so on. The final decision is made on the basis of the results obtained.
  We do not consider this two-step procedure in detail, but note that time-adaptive testing can be applied in this situation, too.

\subsection{The scheme of the time-adaptive testing.} 
Let there be an RNG which  generates binary sequences, and a battery of $s$ tests 
with statistics $\tau_1, \tau_2, ..., \tau_s$.  In addition, suppose that the total available 
testing time is limited to a certain amount $ T $ and the level of significance is  $\alpha \in (0,1)$.

    When the time-adaptive testing is applied, all the calculations are separated into a preliminary stage  and a final one.
The result of the preliminary stage is the list of  values 
$$ 
\gamma_1 = \frac{- \log \pi_{\tau_1}(x_1^1 x_2^1 ... x^1_{n_1})}{n_1} ,  \gamma_2 = \frac{- \log \pi_{\tau_2}(x_1^2   x_2^2 ... x^2_{n_2})}{n_2} 
$$ \begin{equation}\label{prel}
, ... , \, \, \gamma_s = \frac{- \log \pi_{\tau_s}(x_1^s   x_2^s... x^s_{n_s})}{n_s} \,  ,
\end{equation}
where the sequences $ x_1^1 x_2^1 ... x^1_ {n_1} $, $ ... , x_1^s x_2^s ... 
x^s_{n_s}  $  may have common parts (for example, the first sequence may be the prefix of the second, etc.).
Then, taking into account the values (\ref{prel}), it is possible to choose some tests from the battery and apply
them to the longer sequence, calculate new values $\gamma$,  and so on. 
When the preliminary stage is carried out, several tests from the battery should be chosen for the next stage.

The final stage is as follows. First, we  divide    the significance level $\alpha$ into $\alpha_1, \alpha_2, ... , 
\alpha_k$ in such a way that $\sum_{i=1}^k \alpha_i = \alpha$.
Then, we   obtain 
new sequence(s)  $y_1^1 y_2^1 ... y^1_{m_1}, ... , y_1^k   y_2^k  ... y^k_{m_k} $, 
which may have common parts, but are independent of $ x_1^1 x_2^1 ... x^1_ {n_1} $, $ ... , x_1^s x_2^s ... 
x^s_{n_s}  $ and calculate 
\begin{equation}\label{final}
 \pi_{\tau_{i_1}}(y_1^1 y_2^1 ... y^1_{m_1})  
, ... ,  \pi_{\tau_{i_k}}(y_1^k   y_2^k  ... y^k_{m_k}) \,  .
\end{equation}
The hypothesis $H_0$ will be  accepted, if $  \pi_{\tau_{i_j}}(y_1^j y_2^j ... y^j_{m_j})> \alpha_j$ for all $j=1, ... , k$.
Otherwise, $H_0$ is rejected.  The parameters of the test should be chosen in such a way that the total time of calculation is not grater than the given limit $T$.

\begin{claim}\label{t1} 
The  significance level of  the described time-adaptive test is not larger than $\alpha$.
\end{claim}
Indeed, the sequences $y_1^1 y_2^1 ... y^1_{m_1}  
, ... ,  y_1^k   y_2^k  ... y^k_{m_k} \,  $ and $ x_1^1 x_2^1 ... x^1_ {n_1} $, $ ... , x_1^s x_2^s ... 
x^s_{n_s}  $  are independent and, hence,  the results of the final stage does not depend
on the preliminary one.  When the battery $\tau_{i_1}, \tau_{i_2} , ... , \tau_{i_k}$
is applied, the 
significance level of $\tau_{i_j}$ equals $\alpha_{j}$ and the significance level of the battery equals $ \sum_{i=1}^k \alpha_i $. From (\ref{a-g}) we can see that 
the significance level of the battery (and, hence, of the described testing)  is not grater than $\alpha$. 

 {\em Comment.} The length of the sequences may depend on the speed of tests.   
 For example, it can be done as follows: let $v_i $ 
 be the speed per bit of the test $\tau_i$, $i = 1, ... , s$.  One possible way to take into account the speed difference is to calculate
 $$
\hat{ \gamma}_i = \frac{- \log \pi_{\tau_i}(x_1^i x_2^i ... x^i_{n_i})}{n_i / v_i}, \, \, \, i=1, ... , s, 
 $$
 instead of (\ref{prel}) and similar expressions.

\subsection{The experiments.}
We carried out some experiments with the time-adaptive test basing on the   battery
Rabbit from \cite{l2007testu01}, which 
contains 26 tests.  
Let us first describe the choice of the RNG for our experiments.   Nowadays there are many  ``bad'' PRNGs and ``good'' ones. In other words, the output sequences of some known PRNGs have some deviations from randomness, which are quite easy to detect with many known tests, while other PRNGs do not have deviations that can be
detected by known tests \cite{l2007testu01}.  So, 
we need to have some families of RNGs with such deviations from randomness that they can be detected only for quite large output sequences. To do this,  we take a good generator
MRG32k3a and a bad one LCG  from \cite{l2007testu01}, generate sequences 
$g_1 g_2 ... $ and $b_1 b_2 ... $  by these two  generators and then prepared a "mixed" sequence
$m_1 m_2 ... $ in such a way that 

\begin{equation}\label{D}  
 m_i = 
\begin{cases}
    g_i       & \quad \text{if } i \mod D    \, \, \neq 0 \\
   b_i  &    \quad \text{if }     i \mod D    \, \, = 0
    \end{cases}
\end{equation}
where $D$ is a parameter.  

The time-adaptive testing was organised as follows: 
during the preliminary stage we first generated a file $m_1 m_2 ... m_{l_1}$
with $l_1 = 2\, 000\, 000$  bytes, 
   tested it by 
25 tests from the  Rabbit battery and calculated the values (\ref{prel}) with $\log \equiv \log_2$, see the left part of Table 1. (This battery contains 26 tests, but one of them cannot be applied to such a short sequence.)  Then we chose 5 tests with the biggest value
$- log  \, \pi_{t_i}(m_1 ... m_{l_i}) /  l_1 $  (let they be $t_{i_1}, ... , t_{i_5}$ ), generated a sequence
$m_1 ... m_{l_2}$ with $l_2 = 6\, 000 \, 000$  bytes and applied the tests $t_{i_1}, ... , t_{i_5}$ for testing this sequence (see the example in the right part of  Table 1). After that we found a test $t_f$ for which  
$$
- log  \, \pi_{t_f}/  l_f = \max_{r= 1, ..., 25; \, j = i_1 ... i_5}\{ - log \pi_{r}(m_1 ... m_{l_1})/  l_1,  - log \pi_{j}(m_1 ... m_{l_2})/  l_2 \} \, .
$$
 (In other words, 
for $t_f$ the value  $- log \pi_{r}(m_1 ... m_{l_k})/  l_k$ is maximal for $k =1,2$ and
all $r$  (see the  Table 1).
The preliminary stage was finished. Then, during the second stage, we generated a
$ 40\, 000 \,000$ byte  sequence,   and applied the test $t_f$ to it. 
If the obtained p-value was less than $0.001$, the hypothesis $H_0$ was rejected.
(Note that the sequence length  $l_1  = 2\, 000 \, 000$ and  $l_2 = 6\, 000 \, 000$ are 
$5\%$ and $15\%$ from the final length of  $40\, 000 \, 000$ bytes. So, the total 
length of the sequences  tested by all the tests  during the preliminary stage is 
$25 \times 0.05 + 5 \times 0.15$ $= 2$ the final  length, i.e.   $2 \times 40\, 000 \, 000$.
On the other hand, if one applies the battery Rabbit to the sequence of the same length,
 the total length of investigated sequences is $25 \times 40\, 000 \, 000$, i.e. 8,33 times more.

Let us consider one example in detail, taking $D=2$ in (\ref{D}).

\begin{table}[h]
\caption{Time-adaptive testing. Preliminary stage.}
\vspace*{0.4cm}
\begin{tabular}{|p{0.5in}|p{0.5in}|p{0.5in}| p{0.6in}| |p{0.5in}|p{0.5in}|p{0.6in}|} \hline 
test & length ($l$) (bytes) &p-value ($\pi$) & $-\log_2\pi/l$  & length ($l$) (bytes) &  p-value&  $-\log_2\pi/l$
\\ \hline 
 t1 & $2\, 10^6$&0.42 &$6.3\, \, 10^{-7}$	 &   & &   \\ \hline 
 t2 & $2\, 10^6$&0.37 &	$7.3\, \, 10^{-7}$ &  & &   \\ \hline 
 t3 & $2\, 10^6$&0.028 &	$26 \, \, 10^{-7}$ & $6\, 10^6$ &0,23 & 	$3.6 \, \, 10^{-7}$   \\ \hline 
 t4 & $2\, 10^6$& 0.78&	$1.8\, 10^{-7}$ &   & &   \\ \hline 
 t5 & $2\, 10^6$&0.4 &	$6.5\, 10^{-7}$ &   & &   \\ \hline 
 t6 & $2\, 10^6$& 0.37&$7.2 \, 10^{-7}$	 &   & &   \\ \hline 
 t7 & $2\, 10^6$&0.059 &	$20\, 10^{-7}$ &   & &   \\ \hline 
 t8 & $2\, 10^6$&0.026 &	$26\, 10^{-7}$ & $6\, 10^6$ &0.0037 &$26\, \, 10^{-7}$   \\ \hline 
 t9 & $2\, 10^6$&0.72 &	$2.4\, 10^{-7}$ &   & &   \\ \hline 
 t10 & $2\, 10^6$& 0.72&	$2.4\, 10^{-7}$ &   & &   \\ \hline 
 t11 & $2\, 10^6$&0.63 &	$3.3\, 10^{-7}$ &   & &   \\ \hline 
 t12 & $2\, 10^6$&0.74 &	$2.2\, 10^{-7}$ &  & &   \\ \hline 
 t13 & $2\, 10^6$& 0.021&	$28\, 10^{-7}$ & $6\, 10^6$ & 0.0028&  $14\, 10^{-7}$ \\ \hline 
 t14 & $2\, 10^6$&0.42 &	$6.2\, 10^{-7}$  & &  & \\ \hline 
 t15 & $2\, 10^6$& 0.9& $0.74\, 10^{-7}$	 &   &  &  \\ \hline 
 t16 & $2\, 10^6$&0.087 &	$18\, 10^{-7}$ &   & &   \\ \hline 
 t17 & $2\, 10^6$&0.72 &	$2.3\, 10^{-7}$ &   & &   \\ \hline 
 t18 & $2\, 10^6$& 0.58&	$3.9\, 10^{-7}$ &   & &   \\ \hline 
 t19 & $2\, 10^6$&0.89 &	$0.81\, 10^{-7}$ &  & &   \\ \hline 
 t20 & $2\, 10^6$&0.51 &	$4.9\, 10^{-7}$ &   & &   \\ \hline 
 t21& $2\, 10^6$& 0.047&	$22\, 10^{-7}$ & $6\, 10^6$ & 0.73&  $0.76\, 10^{-7}$ \\ \hline 
 t22& $2\, 10^6$& 0.47&	  $0.47\, 10^{-7}$ &  & &   \\ \hline 
 t23& $2\, 10^6$&0.18 &	$12\, 10^{-7}$ &  & &   \\ \hline 
 t24& $2\, 10^6$& 0.14&	$14\, 10^{-7}$ &   & &   \\ \hline 
 t25& $2\, 10^6$& 0.024&	$27\, 10^{-7}$ & $6\, 10^6$ & 0.05& $  7.2\, 10^{-7}$ \\ \hline 
 \end{tabular}
\end{table}

 Table 1 contains the results of all the calculations carried out during the preliminary stage.
So, we can see that the  value $- \log_2 \pi)/l $ is maximal for the test $t13$.
Hence, at the final stage we applied the test $t13$ to the new $40\,  000\, 000$-byte
sequence. It turned out that $\pi_{t13}  =  $ $2.9\, \,10^{- 26}$ and, hence,
$H_0$ is rejected.  Besides, we estimated time of all calculation (during both stages). 

After that, we conducted an additional experiment to get the full picture. Namely, we calculated p-values for all tests and for the same $40\,  000\, 000$-byte
sequence and the estimated the total time of calculations.  It turned out that the p-values of the two tests were less than 0.001. 
Namely,  $\pi_{t13}  =  $ $2.9\, \,10^{- 26}$, $\pi_{t22}  =  $ $1.1\, \,10^{- 6}$.
Besides, we estimated time of calculations for all experiments. 
So, the described time-adaptive testing revealed one of the two most powerful tests,
while the time used is  8 times.

We carried out similar experiments 20 times for $d=  2, 3, 4$ (in (\ref{D}) )  with different good and bad generators from \cite{l2007testu01}.
Besides, we  investigated several modifications
of the considered scheme. In particular, we considered a case where during the  
preliminary  stage we, as before,  first chose 5 the best tests and them two of the best  tests for the finale stage (instead of one, as in the  experiment above).
 It turned out, that in all  cases  considered the battery Rabbit rejects $H_0$  and  the time-adaptive testing rejected  $H_0$,  too.

\section{Conclusion}
First of all, we note that the proposed time-adaptive testing does not suggest exact values of numerous parameters. Among these parameters, we note the number of steps at the preliminary stage (in the considered example there were two such steps: selecting five tests and then one), the number of tests compared in one step, the length of the tested sequences, the rule for choosing tests at different stages, etc.
The problem of choosing the parameters may be considered a problem of multidimensional optimization. There are many methods available for solving such problems (for example, neural networks and other AI algorithms), and some of them can be used along with the time-adaptive testing.

As far as we know, no one has applied adaptive methods for testing randomness, but there are several well-known approaches that can be considered as steps in this direction. For example, L'Ecuyer and Simard recommend several batteries of different sizes that require different times (and the investigator may use them depending on how much time he has) \cite{l2007testu01}. Another popular battery recommended for cryptographic applications also has some parameters that allow one to adjust the testing time \cite{NIST-test}.

We believe that the proposed approach makes it possible to investigate and optimize time-adaptive testing. 
\section{Appendix 1. Consistent  tests based on universal codes.}

The considered  tests are based on so-called universal codes, that is why we first briefly 
describe them. 
For any integer $m$ a  code $\phi$ is defined  as such a map  from the set
of $m$-letter words to the set of all binary words that for any 
$m$-letter $u$ and $v$ $\, \, \,  \phi(u) \neq \phi(u)$. This property gives a possibility
to uniquely decode. 
(More formally, $\phi$ is injective mapping from 
$\{0,1\}^m$ to $\{0,1\}^*$, where  $\{0,1\}^* = \bigcup _{i=1}^\infty	\{0,1\}^i$.)
We will consider so-called universal codes which have the two  following properties:
 \begin{equation}\label{un-krafr}
\forall \; m >0  \quad  \sum_{u \in \{0,1\}^m} 2^{- | \phi(u) | } =1\,
\end{equation} 
and for any stationary ergodic $\nu$  defined on the set of all infinite binary words
$x   = x_1 x_2 ... $, with probability one
 \begin{equation}\label{un-lim}
\lim_{n \to \infty } \frac{1}{n} | \phi( x_1 x_2 ... x_n) | / n  = h(\nu) \,
\end{equation}
where $ h(\nu) $ is the Shannon entropy of $\nu$. Such code exist,  see \cite{co}.
Note, that a goal of  codes is to " compress " sequences, i.e. make an average length of the codeword 
$\phi( x_1 x_2 ... x_n) $ as small as possible. The second property (\ref{un-lim}) shows
that the universal codes are asymptotically optimal, because the Shannon entropy is 
a low bound of the length of the compressed sequence (per letter), see \cite{co}.

Let us back to considered problem of hypothesis testing. 
Suppose, it is known that a sample sequence $x= x_1 x_2 ...$  was  generated by stationary ergodic source and, as before,  we consider the same $H_0$ against the same $H_1$. Let $\phi$ be a universal code. 
The following test is suggested in \cite{R1}:

If the length  $|\phi (x_1 ... x_n)$  $\le n - \log_2 \alpha $ then $H_0$ is rejected, 
otherwise accepted.
Here, as before, $\alpha$ is the significance level,  $|\phi (x_1 ... x_n)|$  is the length 
of encoded (''compressed") sequence. We denote  this test by $T_\phi$ and its statistic  by $\tau_\phi$, i.e. 
 \begin{equation}\label{tau-fi}
\tau_\phi (x_1 ... x_n) = n - |\phi(x_1 ... x_n)| \, .  
\end{equation}

The following theorem is proven in \cite{R1,R2}:
\begin{theorem}\label{t0}  
 For each stationary ergodic $\nu$, $alpha \in (0,1) $ and a universal code $\phi$, with probability 1 the Type I error of the described test is not larger than  $\alpha$ and  the Type II error goes to 0, when $n \to \infty$.
\end{theorem}

\section{Appendix 2.  Proofs.}
{\it Proof of Theorem 1.}   
 The known  Shannon-McMillan-Breiman (SMB) theorem   claims that for the stationary ergodic source
 $\nu$ and any $\epsilon >0, \delta > 0$ there exists such $n'$ that 
 \begin{equation}\label{smb}
 \nu \{ x: x\in \{0,1\}^n  \,\, \,  \& \, \,  \, h(\nu) - \epsilon < - \frac{1}{n } \log \nu(x)   <
 h(\nu) + \epsilon \, \, \} >1 - \delta \, 
 \end{equation} 
 for $n > n'$, 
 see \cite{co}.  From this we obtain 
 \begin{equation}\label{smb2}
 \nu \{  x: x\in \{0,1\}^n  \,\, \,  \& \, \,   2^{- n (h(\nu) -  \epsilon )} >\, \nu(x) >  2^{- n (h(\nu) +  \epsilon )}  \}
 >1 - \delta \, 
 \end{equation}
  for $n > n'$. 
  It will be convenient to define
  \begin{equation}\label{fi}
\Phi_{\epsilon, \delta, n} = \{ x: x\in \{0,1\}^n  \,\, \,  \& \, \,  \, h(\nu) - \epsilon < - \frac{1}{n } \log \nu(x)   <
 h(\nu) + \epsilon \, \, \}  
 \end{equation} 
  From this definition and (\ref{smb2} )  we obtain 
 \begin{equation}\label{f2}
(1-\delta) \, 2^{n (h(\nu) -\epsilon) } \le | \Phi_{\epsilon, \delta, n} | \le 2^{n (h(\nu) +\epsilon) }\, .
 \end{equation} 
For any $x \in \Phi_{\epsilon, \delta, n}$ define 
 \begin{equation}\label{lam}
\Lambda_x  = \{ y: \nu (y) > \nu(x) \, \, \}  \bigcap  \Phi_{\epsilon, \delta, n} \, .
 \end{equation}  Note that, by definition, $|\Lambda_x| \le | \Phi_{\epsilon, \delta, n}|$
 and from (\ref{f2}) we obtain
 \begin{equation}\label{lam2} | \Lambda_x | \le 2^{n (h(\nu) +\epsilon) }\, . \end{equation} 
For any $\rho \in (0,1)$ we define $\Psi_\rho \subset \Phi_{\epsilon, \delta, n}$ such that
 \begin{equation}\label{psi}
 \nu( \Psi_\rho) = \rho \,\, 
\&  \,  \, \forall 	u \in \Psi_\rho \,    \forall 	v \in ( \Phi_{\epsilon, \delta, n} \setminus 
\Psi_\rho )\,
\to
 \nu(u) )\ge \nu(v) \, .
 \end{equation} 
 (That is, $\Psi_\rho $ contains the most probable words whose total probability equals $\rho$.)
 Let us consider any $x \in ( \Phi_{\epsilon, \delta, n} \setminus 
\Psi_\rho )  $.  Taking into account 
the definition (\ref{psi}) and (\ref{f2})
we can see that for this $x$
  \begin{equation}\label{lval-}
|\Lambda_x | \ge \rho | \Phi_{\epsilon, \delta, n}| \ge \rho (1-\delta) 2^{n (h(\nu) -\epsilon) } \, .
 \end{equation}  
So, from this inequality and (\ref{lam2}) we obtain 
\begin{equation}\label{lam-size}
\rho (1-\delta) 2^{n (h(\nu) -\epsilon) } \le |\Lambda_x | \le \,2^{n (h(\nu) +\epsilon) } \, . \end{equation}   
From equation  (\ref{smb2}), (\ref{fi}) and (\ref{psi}) we can see that 
$\nu ( \Phi_{\epsilon, \delta, n} \setminus \Psi_\rho ) \ge (1-\delta) (1-\rho) $.
Taking into account (\ref{lam-size}) and this inequality, we can see that 
$$
 \nu \{x: x \in \{0,1\}^n \& 
 h(\nu) - \epsilon - \log ( \rho (1-\delta))/ n 
\le  \log
|\Lambda_x | /n  
 $$
  \begin{equation}\label{fin}
  \le h(\nu) +  \epsilon \} \ge  (1-\delta) (1-\rho) . 
 \end{equation}  
 From the definition (\ref{pi-np}) of $\pi_{NP}(x)$ and the definition
 (\ref{lam}) of $\Lambda_x $, we 
can see that $ \pi_{NP}(x) = |\Lambda_x |/ 2^n$.
Taking into account this equation and (\ref{fin})
 we obtain  the following:
 $$  \nu \{x: x \in \{0,1\}^n \, \& \,1 - ( h(\nu) - \epsilon - \log  ( \rho (1-\delta)) / n ) \ge $$
 \begin{equation}\label{fin2} 
  - \log
\pi_{NP}(x) /n  \ge 1 - ( h(\nu) +  \epsilon )  \} \ge  (1-\delta) (1-\rho) . 
 \end{equation}  
 Having taken into account that this inequality is valid for all positive $\epsilon, 
 \delta$ and $\rho$ we obtain the first statement of the theorem.

The proof of the second statement of the theorem  is closed to the previous one.  First,
from the theorem 2 we  see that  for any
$\epsilon >0, \delta >0$ 
 we define 
 \begin{equation}\label{abc}
 \hat{\Phi}_{\epsilon, \delta, n} = \{x: h(\nu) - \epsilon < | \phi(x_1 ... x_n) | /n <  h(\nu) + \epsilon \, \}\, .
 \end{equation}  
Note that
from (\ref{un-lim} ) we can see that there exists  such $n''$ that,
for $n > n''$, 
\begin{equation}\label{fi-hat}
\nu ( \hat{\Phi}_{\epsilon, \delta, n})  > 1 - \delta \, .
 \end{equation}  
 We will use the set $\Phi_{\epsilon, \delta, n}$ (see (\ref{fi}) ). 
 Having taken into account 
the  SMB theorem (\ref{smb}) and   (\ref{fi-hat}), 
we can see that 
\begin{equation}\label{fi-hat2}
\nu ( \hat{\Phi}_{\epsilon, \delta, n} \cap    \Phi_{\epsilon, \delta, n} )    > 1 - 2 \delta \, ,
 \end{equation}  
if $n > \max (n', n'')$. 

From this moment, the proof begins to repeat the proof of the first statement if we use the set $(\hat{\Phi}_{\epsilon, \delta, n} \cap    \Phi_{\epsilon, \delta, n} ) $ instead of  $ \Phi_{\epsilon, \delta, n}$.  The only difference is in the definitions (\ref{lam}) and (\ref{psi}) which should be changed as follows.
$$
\Lambda_x  = \{ y: |\phi| (y) | < |\phi(x)| \, \, \}  \cap  (\hat{\Phi}_{\epsilon, \delta, n} \cap    \Phi_{\epsilon, \delta, n} ) \, 
$$
and $ \Psi_\rho $ is such a subset of 
$(\hat{\Phi}_{\epsilon, \delta, n} \cap    \Phi_{\epsilon, \delta, n} )$  that
$$
 \nu( \Psi_\rho) = \rho \,\, 
\&  \,  \, \forall 	u \in \Psi_\rho \,    \forall 	v \in ( \Phi_{\epsilon, \delta, n} \setminus 
\Psi_\rho )\, 
\to
 |\phi(u) |\le |\phi(v) | \, .
$$
If we replace $\pi_{NP}$ with  $\pi_{\tau_\phi}$ and $\delta$ with $2 \delta$, we
obtain the proof of the second statement.
Theorem is proven.

\section*{Acknowledgment}
This work was supported by Russian Foundation for Basic Research (grant 18-29-03005).

\end{document}